\documentclass[10pt]{article}

\pdfoutput=1

\usepackage{amssymb}
\usepackage{amsmath}
\usepackage{amsthm}
\usepackage{graphicx}
\usepackage{subfigure}
\usepackage{dsfont}
\usepackage{mathrsfs}
\usepackage[all]{xy}
\usepackage{rotating}
\usepackage{hyperref}

\newcommand{\bft}{{\bf t}}

\newcommand{\mapJ}{\mathscr{J}}

\renewcommand{\labelenumi}{(\alph{enumi})}

\newenvironment{packed_enum}{
\begin{enumerate}
  \setlength{\itemsep}{1pt}
  \setlength{\parskip}{0pt}
  \setlength{\parsep}{0pt}
}{\end{enumerate}}

\theoremstyle{definition}
\newtheorem{theorem}{Theorem}

\newtheorem{definition}[theorem]{Definition}

\newtheorem*{theorem*}{Theorem}
\newtheorem*{corollary*}{Corollary}
\newtheorem*{definition*}{Definition}
\newtheorem*{example*}{Example}
\newtheorem*{exercise*}{Exercise}
\newtheorem*{lemma*}{Lemma}
\newtheorem*{proposition*}{Proposition}
\newtheorem*{remark*}{Remark}
\newtheorem*{conjecture*}{Conjecture}

\title{Attractive $n$-Type Contact Processes}
\author{Joseph Stover\footnote{jstover@bren.ucsb.edu}}
\begin{document}
\maketitle
\let\thefootnote\relax\footnotetext{This work was partially supported by the NSF Vigre graduate fellowship, and the CATTS fellowships at the University of Arizona}
\let\thefootnote\relax
\let\thefootnote\relax
%%%%%%%%%%%

\begin{abstract}
Interacting particle systems are continuous time Markov processes which are used to construct stochastic spatial models. Monotonicity is a useful property which simplifies certain calculations, one of which is the ability to use computational algorithms to sample exactly from the stationary distribution for certain processes. A monotone interacting particle system is called attractive. Monotonicity is well understood for spin systems which only include two particle types, such as the contact process, however, when constructing applied models, it is often desirable to include more. In this paper, an interaction map is used to describe the interactions that occur in a model and to understand monotonicity for a certain class of multitype contact processes.
\end{abstract}

%\begin{keyword}[class=AMS]
%\kwd{60K35}
%\end{keyword}
%
%\begin{keyword}
%\kwd{monotone}
%\kwd{attractive}
%%\kwd{spin system}
%\kwd{coupling}
%\kwd{multitype contact process}
%\kwd{interacting particle system}
%\end{keyword}

%\end{frontmatter}

\section{Introduction}
Interacting particle systems are a class of Markov processes used to model the evolution of particles types $S$ on a collection of sites $\Lambda\subseteq\mathbb{Z}^d$. Consequently, a state for these processes is a mapping $\eta: \Lambda \rightarrow S$ that
assigns a particle type to each site. A process realization $\eta_{\bf t}=\{\eta_t: t\in[0,\infty)\}$ is a right continuous function of $t$. Elements of the state space, $\Xi=S^\Lambda$, are interchangeably referred to as states or configurations. If $\eta(x)=a$, site $x$ is said to be infected or occupied by particle type $a$, or more simply, particle $a$ sits at site $x$. Typically $0 \in S$ and $\eta(x) = 0$ indicates an unoccupied site.

A general interacting particle system is quite difficult to analyze and some additional properties must be considered in order to make meaningful general statements. In this paper, our motivating examples are multi-species biological models. Thus, we will be primarily focused on circumstances in which the number of particles types $\#|S|$ is at least three -- vacant and two species of organism. %The sites $\Lambda$ in this example are finite subsets of $\mathbb{Z}^d$.

The contact process \cite{Harris74} is a basic interacting particle system used to build many biological models, and the nature of the interactions which occur has led us to introduce the `interaction map' which characterizes the allowable interactions in a model.

The main result is a characterization of those interaction maps that lead to monotone processes. In this case, the particle types $S$ must be ordered thus endowing the state space with a partial ordering `$\le$' (sitewise comparison). Monotonicity allows us to reduce the analysis of many properties of an interacting particle system to a study of the process considering the evolution starting from the small number on initial configurations that are extremal under the given partial order. For example, the coupling from the past (CFTP) algorithm \cite{Propp96} allows us to obtain an exact sample from the stationary distribution of an irreducible and aperiodic Markov chain with a finite state space. For monotone processes, the algorithm need only be applied to the extremal states $\bar{0}$ and $\bar{n}$ (the states where all sites are identically occupied by $0$ or $n$ particles respectively) such that $\bar{0}\le\eta\le\bar{n}$ for all $\eta\in\Xi$.

The interaction map formulation also allows for a fast assessment of whether or not reordering the particles may make the process attractive. We show that switching the order of the basic two type contact process results in attractiveness.

We will briefly review the contact process and monotonicity. This will set the stage for the development of the interaction map.

\section{Contact Processes}
%The rate at which a particle at a site on the lattice is replaced by another particle type depends on the types of particles occupying the nearest neighboring sites. Different types of neighborhood structures are often used, and the results presented here generalize to any size or shape of neighborhood.

The contact process is a {\em spin system} (it only has two particle types), and it is the basic interacting particle system of interest here. Its state space is $\{0,1\}^{\mathbb{Z}^d}$ with particle type 0 representing an empty site and type 1 an occupied site. The transition rates are:
\begin{equation*}
\begin{tabular}{ccc}
transition & rate & description\\
$0\rightarrow1$ & $\lambda n_1$ & birth\\
$1\rightarrow0$ & $1$ & death\\
\end{tabular}
\end{equation*}
where $n_1$ is the number of occupied neighbors. A birth  occurs at a rate proportional to the number of infected neighbors. An occupied site becomes vacant at rate one. The transition rate for a spin system is referred to as the {\em flip rate}, and denoted $c(x,\eta)$. Thus for the contact process,
\begin{equation}
c(x,\eta) =\left\{
\begin{aligned}
&\sum_{y:|x-y|=1} I_{1}(\eta(y)) &\text{ if } \eta(x)=0 \\
&1 &\text{ if } \eta(x)=1
\end{aligned}\right.
\end{equation}
where $I_{1}(\cdot)$ is the indicator function for particle type 1.

Numerous variations have been created around the rule that infection occurs at a rate proportional to the number of occupied neighboring sites. A contact process which has $n+1$ particle types will be referred to as an {\em $n$-type contact process} (usually type 0 denotes an unoccupied site). The basic two-type contact process is a competition model studied by Neuhauser \cite{NeuDiss}. This process has the following transition rules:
\begin{equation*}
\begin{tabular}{cccc}
$0\rightarrow1$ & $\lambda_1 n_1$ & $0\rightarrow2$ & $\lambda_2 n_2$\\
$1\rightarrow0$ & $\delta_1$ & $2\rightarrow0$ & $\delta_2$\\
\end{tabular}
\end{equation*}
The two-stage contact process studied by Krone \cite{mScp_Krone} is another two-type contact process. It is a single species model which includes two life stages. The transitions are such that 2's give birth to 1's and the 1's mature into 2's at a constant rate:
\begin{equation*}
\begin{tabular}{cccc}
$0\rightarrow1$ & $\lambda n_2$ & $1\rightarrow2$ & $\gamma$\\
$1\rightarrow0$ & $1+ \delta$ & $2\rightarrow0$ & $1$\\
\end{tabular}
\end{equation*}
The grass--bushes--trees successional model proposed by Durrett and Swindle \cite{DurSwin91} is a two-type contact process as well. It is the basic two-type contact process with the modification that 2's are allow to give birth onto sites occupied by 1's in addition to empty sites:
\begin{equation*}
\begin{tabular}{cccccc}
$0\rightarrow1$ & $\lambda_1 n_1$ & $0\rightarrow2$ & $\lambda_2 n_2$ & $1\rightarrow2$ & $\lambda_2 n_2$\\
$1\rightarrow0$ & $\delta_1$ & $2\rightarrow0$ & $\delta_2$\\
\end{tabular}
\end{equation*}

These models are all multitype contact processes, however they are no longer spin systems so what is generally required for attractiveness has not been previously known. In \cite{mScp_Krone} (Theorem 3.1), it is shown that the two-stage contact process is monotone with respect to its parameters. In \cite{CoexComp_DurNeu} (Proposition 1.1), a more complicated type of monotonicity property is described for a three-type contact process.

\section{Monotonicity}
A Markov process with a partially ordered state space, $\Xi$, and semigroup, $S(t)$, is called {\em monotone} if either of the equivalent conditions, (\ref{eq:monSemi}) or (\ref{eq:monCoup}), is shown to be satisfied.
\begin{subequations}
\begin{eqnarray}
f\in\mathscr{M} \text{ implies } S(t)f\in\mathscr{M} \text{ for all } t\ge0 \label{eq:monSemi}\\
\mu_1\le\mu_2 \text{ implies } \mu_1S(t)\le\mu_2S(t) \text{ for all } t\ge0 \label{eq:monCoup}
\end{eqnarray}
\end{subequations}
The $\mu_i$ are probability measures on the state space, and $\mathscr{M}$ is the set of continuous monotone functions, $f:\Xi\rightarrow\mathbb{R}$, such that two states satisfying $\eta\le\xi$ implies $f(\eta)\le f(\xi)$. For the proof of equivalency, see \cite{LiggettIPS} (Chapter 2, Theorem 2.2). Furthermore, $\mu_1\le\mu_2$ is equivalent to there existing a measure, $\nu$ on $\Xi\times\Xi$, that satisfies
\begin{packed_enum}
\item $\nu\{(\eta,\xi) : \eta\in A\}=\mu_1(A)$, and
\item $\nu\{(\eta,\xi) : \xi\in A\}=\mu_2(A)$, where $A$ is any Borel set in $\Xi$, and
\item $\nu\{(\eta,\xi) : \eta\le\xi\}=1$.
\end{packed_enum}
%In other words, the ability to construct a coupling that preserves the partial ordering a.s. is equivalent to the process being monotone.
For the proof of this, see \cite{LiggettIPS} (Chapter 2, Theorem 2.4). %need corollary 1.7,1.8 from liggett.
For interacting particle systems, the typical route is to define what it means for the process to be attractive and then to prove that this is equivalent to being monotone.

%\subsection{Attractive Spin Systems}
For a spin system, given $\eta\le\xi$, attractiveness is equivalent the the flip rate satisfying
\begin{subequations}
\begin{eqnarray}
& c(x,\eta) \le c(x,\xi) &\text{when } \ \eta(x)=\xi(x)=0, \label{eq:spinupRateRest}\\
& c(x,\eta) \ge c(x,\xi) &\text{when } \ \eta(x)=\xi(x)=1. \label{eq:spindownRateRest}
\end{eqnarray}
\end{subequations}
%This means that given any state, replacing some of the 0's in the configuration of particles by 1's can only maintain or increase the flip rate of sites occupied by 0's, and should maintain or decrease the flip rate for sites occupied by 1's.
See \cite{LiggettIPS} (Chapter 3, Theorem 2.2) for the short proof that (\ref{eq:spinupRateRest}) and (\ref{eq:spindownRateRest}) are together equivalent to monotonicity. This paper generalizes these conditions to a class of multitype contact processes.

\subsection{The Interaction Map}
The approach here is to define what is called the {\em interaction map}. We start with a finite set of totally ordered particle types, $S=\{0,1,2,\ldots,n\}$ and define a mapping that describes all interactions between them.
\begin{definition} Given a finite set of totally ordered particle types, $S$, the map, $\mathscr{J} : S\times S\rightarrow S$ is called an {\em interaction map} if its domain is all of $S\times S$, and its range is a subset of $S$. The particle type that replaces type $a$ upon interaction with type $b$ is given by $\mathscr{J}(a,b)$, and this interaction is denoted by $a\overset{b}{\rightarrow}\mathscr{J}(a,b)$.
\end{definition}
The domain of $\mathscr{J}$ is partitioned into
sets %three {\em interaction sets}
of {\em up}, {\em null}, and {\em down} interactions, respectively: $\mathcal{U}=\{(a,b) \in S\times S | \mathscr{J}(a,b) > a\}$, $\mathcal{N}=\{(a,b) \in S\times S | \mathscr{J}(a,b) = a\}$, and $\mathcal{D}=\{(a,b) \in S\times S | \mathscr{J}(a,b) < a\}$. %This partitions the interactions into those which result in larger, unchanged, and smaller particle values respectively.
A process whose interactions are completely defined by a single interaction map is referred to as a {\em interaction map system} (IMS). Many processes may also be described using multiple interaction maps; this is discussed briefly in the last section but is not pursued in detail here.

For the contact process, $\mapJ(0,1)=1$ denotes a birth onto an unoccupied site ($0\overset{1}{\rightarrow}1$), and $\mapJ(1,i)=0$ for any $i$ denotes a death ($1\overset{i}{\rightarrow}0$). The interaction map for the voter model is the same as for the contact process except that $\mapJ(1,1)=1$ because a voter only changes opinions by contact with its opposite (Figure \ref{fig:spinIntMaps}). %($1\overset{\hspace{2px}\text{\scalebox{2}{${}_{_.}$}}}{\rightarrow}$0)
%
%%%%%%%%%%%%%%%%%%%%%%%%%%%%%
% Spin Int Map picture:
\begin{figure}[h]
\begin{center}
\begin{picture}(120,60)
%contact
\put(0,10){\begin{tabular}{| c || c | c | }
\hline
$_b\hspace{-4pt}\rotatebox{135}{\rule{10pt}{0.5pt}}\hspace{-4pt} ^a$  & 0 & 1 \\ \hline\hline
0 & 0 & 0 \\ \hline
1 & 1 & 0 \\ \hline
\end{tabular}}
\put(0,46){\small{particle affected by the neighbor}}
\put(-18,-4){\rotatebox{90}{\small{neighbor}}}
%\put(-18,-4){$b$}
\put(20,-26){(c)}
%voter
\put(80,10){\begin{tabular}{| c || c | c | }
\hline
$_b\hspace{-4pt}\rotatebox{135}{\rule{10pt}{0.5pt}}\hspace{-4pt} ^a$  & 0 & 1 \\ \hline\hline
0 & 0 & 0 \\ \hline
1 & 1 & 1 \\ \hline
\end{tabular}}
\put(100,-26){(v)}
\end{picture}
\end{center}
\vspace{16pt}
\caption{Interaction maps for the (c) Contact Process and (v) Voter model.} %(in table format)}
\label{fig:spinIntMaps}
\end{figure}
% end spin Int Map picture
%%%%%%%%%%%%%%%%%%%%%%%%%%%%%%%
%
The interaction map for the basic two-type contact process is built from that for the contact process. It additionally has $\mapJ(0,2)=2$ and $\mapJ(2,i)=0$ (for any $i$) for births and deaths of species two respectively (Figure \ref{fig:multiCPMap}).
%
%%%%%%%%%%%%%%%%%%%%%%%%%%%%%
% Multi-type CP Map picture:
\begin{figure}[h]
\begin{center}
\begin{picture}(40,40)
%multi-contact
\put(-20,0){\begin{tabular}{| c || c | c | c |}
\hline
$_b\hspace{-4pt}\rotatebox{135}{\rule{10pt}{0.5pt}}\hspace{-4pt} ^a$  & 0 & 1 & 2\\ \hline\hline
0 & 0 & 0 & 0\\ \hline
1 & 1 & 0 & 0\\ \hline
2 & 2 & 0 & 0\\ \hline
\end{tabular}}
\end{picture}
\end{center}
\vspace{14pt}
\caption{Interaction map formulation for the multitype contact process}
\label{fig:multiCPMap}
\end{figure}
% end Multi-type CP Map picture
%%%%%%%%%%%%%%%%%%%%%%%%%%%%%%%

%\begin{definition}
The interaction map, $\mathscr{J}$, is called {\em non-decreasing} on $A\subset S\times S$ if for any $(a_1,b_1)$ and $(a_2,b_2)$ in $A$ which satisfy $(a_1,b_1)\le(a_2,b_2)$ (meaning $a_1\le a_2$ and $b_1\le b_2$), it follows that $\mathscr{J}(a_1,b_1)\le\mathscr{J}(a_2,b_2)$.
%\end{definition}
\begin{definition}
The interaction map will be called {\em attractive} if it satisfies the following conditions:
\begin{packed_enum}
\item $\mathscr{J} \text{ is non-decreasing on } \mathcal{U}$.
\item $\mathscr{J} \text{ is non-decreasing on } \mathcal{D}$.
\item If $(a_1,b_1)\le(a_2,b_2)$, $(a_1,b_1)\in\mathcal{U}$, and $(a_2,b_2) \in \mathcal{D}\cup\mathcal{N}$, then $\mathscr{J}(a_1,b_1)\le a_2$.
\item If $(a_1,b_1)\le(a_2,b_2)$, $(a_1,b_1)\in\mathcal{U}\cup\mathcal{N}$, and $(a_2,b_2) \in \mathcal{D}$, then $a_1\le \mathscr{J}(a_2,b_2)$.
\end{packed_enum}
\label{def:attrMap}
\end{definition}
The last two conditions of Definition \ref{def:attrMap} are equivalent to:
\begin{packed_enum}\setcounter{enumi}{2}\renewcommand{\labelenumi}{{\em(\alph{enumi}\hspace{1px}$^{\prime}$)}}
\item If $(a_1,b_1)\in\mathcal{U}$, then $\mathscr{J}(a_1,b_1)\le \displaystyle\min_a\{a :(a,b)\in \mathcal{D}\cup\mathcal{N}\}$.
\item If $(a_2,b_2)\in\mathcal{D}$, then $\displaystyle\max_a\{a :(a,b)\in \mathcal{U}\cup\mathcal{N}\}\le \mathscr{J}(a_2,b_2)$.
\end{packed_enum}

Definition \ref{def:attrMap} is designed so that given two ordered pairs of particle types, $(a_1,b_1)\le(a_2,b_2)$, we have $\mapJ(a_1,b_1)\le\mapJ(a_2,b_2)$ when both pairs have an up or both a down interaction (or one of them is null)%(Definition \ref{def:attrMap}, (a) and (b))
. When one interaction is up and the other is down%(or one of them is null)
, the ordering need not be preserved, but cannot be broken arbitrarily; it must obey (c) and (d) of Definition \ref{def:attrMap} (or equivalently (c') and (d')). %that if two ordered configurations undergo simultaneous interactions at a site, then the ordering will be preserved.
Staring from Definition \ref{def:attrMap}, a little extra work shows that an attractive interaction map $\mapJ(a,b)$ is nondecreasing in $b$ and is nondecreasing in $a$ except for particle pairs satisfying $(a,b)\in\mathcal{U}$ and $(a+1,b)\in\mathcal{D}$, in which case $\mapJ(a,b)=a+1$ and $\mapJ(a+1,b)=a$.

The interaction map of the basic two-type contact process (Figure \ref{fig:multiCPMap}) is not attractive because $\mathscr{J}(0,2)=2>\displaystyle\min_a\{a:(a,2)\in\mathcal{N}\cup\mathcal{D}\}=1$. There are several ways to modify this interaction map into one which is attractive. Keeping the birth onto empty sites unchanged, necessitates that (a) $\mathscr{J}(1,2)=1$ or (b) $\mathscr{J}(1,2)=2$ and $\mathscr{J}(2,2)\in\{1,2\}$ (Figure \ref{fig:multiCPattrMaps}). Thus the death rates are no longer both constant, or the two species competition character of the model may be broken.
%%%%%%%%%%%%%%%%%%%%%%%%%%%%%
% Multi-type CP Map picture:
\begin{figure}[h]
\begin{center}
\begin{picture}(120,60)
%multi-contact
\put(-80,20){\begin{tabular}{| c || c | c | c |}
\hline
$_b\hspace{-4pt}\rotatebox{135}{\rule{10pt}{0.5pt}}\hspace{-4pt} ^a$  & 0 & 1 & 2\\ \hline\hline
0 & 0 & 0 & 0\\ \hline
1 & 1 & 0 & 0\\ \hline
2 & 2 & 1 & 1\\ \hline
\end{tabular}}
\put(-48,-26){(a)}%when 2 is crowded.}
\put(20,20){\begin{tabular}{| c || c | c | c |}
\hline
$_b\hspace{-4pt}\rotatebox{135}{\rule{10pt}{0.5pt}}\hspace{-4pt} ^a$  & 0 & 1 & 2\\ \hline\hline
0 & 0 & 0 & 0\\ \hline
1 & 1 & 0 & 0\\ \hline
2 & 2 & 2 & 1\\ \hline
\end{tabular}}
\put(56,-26){(b)}%2--2 interaction is null.}%Species 2 immune to crowding.}
\put(120,20){\begin{tabular}{| c || c | c | c |}
\hline
$_b\hspace{-4pt}\rotatebox{135}{\rule{10pt}{0.5pt}}\hspace{-4pt} ^a$  & 0 & 1 & 2\\ \hline\hline
0 & 0 & 0 & 0\\ \hline
1 & 1 & 0 & 0\\ \hline
2 & 2 & 2 & 2\\ \hline
\end{tabular}}
\put(154,-26){(c)}
\end{picture}
\end{center}
\vspace{16pt}
\caption{Attractive modifications to the basic two-type contact process interaction map. The interaction between (1,2) and (2,2) must be modified so that Definition \ref{def:attrMap} is satisfied. The death rates are no longer constant.}
\label{fig:multiCPattrMaps}
\end{figure}
% end Multi-type CP Attractive Map pictures
%%%%%%%%%%%%%%%%%%%%%%%%%%%%%%%

\begin{figure}[h]
\begin{center}
\includegraphics[scale=0.16]{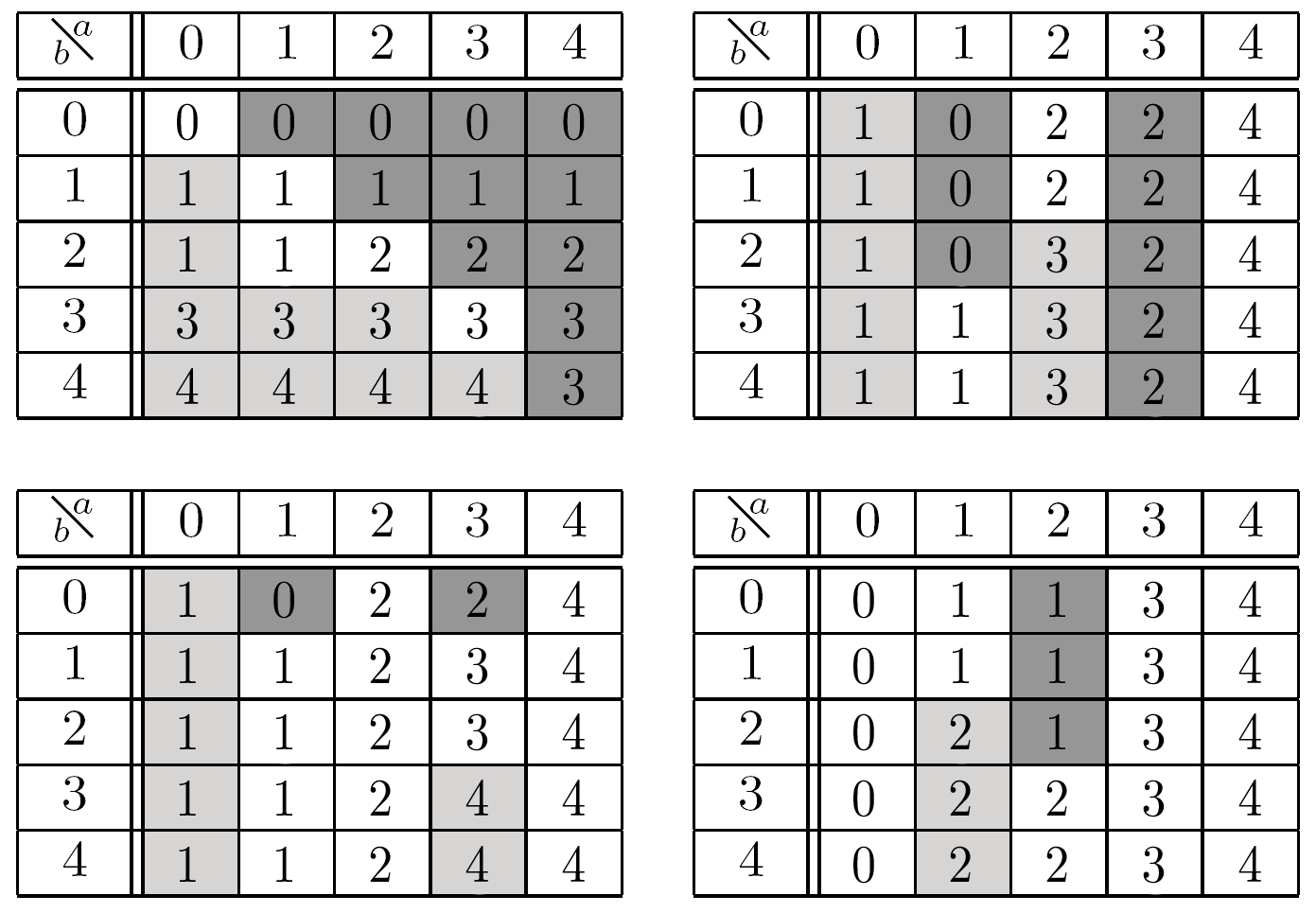}
\end{center}
\caption{Examples of attractive interaction maps are shown above. Those cells shaded dark gray represent down interactions, light gray up interactions, and unshaded cells are null interactions. For example, in the upper left interaction map, $\mathscr{J}(2,3)=3$. This shows that particle type 3 influences type 2 to become type 3 and is thus an up interaction.}
\label{fig:attrMaps}
\end{figure}

Figure \ref{fig:attrMaps} shows four examples of attractive interaction maps with the different set of interactions shaded. Dark gray denotes down interactions, light gray denotes up interactions, and no shading denotes null interactions. We can see that the interaction maps are nondecreasing separately in the light gray and dark gray regions and that when a light gray cell is on the left of a dark gray cell, the decrease is exactly one: $\mapJ(a,b)=a+1$ and $\mapJ(a+1,b)=a$.

%From the example interaction maps in Figure \ref{fig:attrMaps}, it is clear that when an up cell is directly left of a down cell, the interactions are defined by $\mapJ(a,b)=a+1$ and $\mapJ(a+1,b)=a$ for some $b$. The maps are all clearly non-decreasing on the sets $\tU\cup\tN$ and $\tD\cup\tN$. These are the conditions in Theorem \ref{thm:attrMap2}. All cells above a down interaction cell must also represent down interactions, and all cells below an up interaction cell must also represent up interactions(Lemma \ref{lem:attrMap1}). Theorem \ref{thm:stepMapAttr} is illustrated by the fact that each attractive interaction map is stepwise non-decreasing vertically and horizontally except for horizontal decreases by exactly one in accordance with Theorem \ref{thm:attrMap2}, Corollary \ref{cor:attrMap4}, and Theorem \ref{thm:stepMapAttr}.

\section{Transition Rates}
The transition rates are denoted by $r_u(\eta,x,y)$ and $r_d(\eta,x,y)$, which give the rates for up and down interactions respectively. By convention, null interactions are assigned a transition rate of zero since they do not change the configuration. For each pair of particle types $(a,b)\in S\times S$, $\lambda_{ab}$ denotes the rate at which the interaction  $a\overset{b}{\rightarrow}\mathscr{J}(a,b)$ occurs. For a configuration $\eta$ the rates at which a particle at site $y$ influences a particle at site $x$ to change are thus given by
\begin{subequations}
\begin{eqnarray*}
r_u(\eta,x,y)&=&\sum_{a\in S,b\in \mathcal{U}_a}\lambda_{ab}\phi(x,y)I_{\{(a,b)\}}(\eta(x),\eta(y))\label{eq:r_u}\\
r_d(\eta,x,y)&=&\sum_{a\in S,b\in \mathcal{D}_a}\lambda_{ab}\phi(x,y)I_{\{(a,b)\}}(\eta(x),\eta(y))\label{eq:r_d}.
\end{eqnarray*}
\end{subequations}
The function $\phi(x,y)$ defines the neighborhood, and is assumed to be non-negative and bounded over a finite neighborhood. The sets $\mathcal{U}_a=\{b\in S : \mathscr{J}(a,b) > a\}$ and $\mathcal{D}_a=\{b\in S : \mathscr{J}(a,b) < a\}$ for a given particle type $a$ denote the sets of particle types with up or down interactions respectively.

%The parameter $\lambda_{ab}$, is the rate parameter controlling the interaction that particle type $a \overset{b}{\longrightarrow} \mathscr{J}(a,b)$ due to the presence of a particle of type $b$ in its neighborhood. If $\lambda_{ab}\equiv0$ for some $a$ and $b$, then by convention, we put $(a,b)\in \mathcal{N}$, choosing $\mathscr{J}(a,b) = a$, for the process. This is so that all interactions will either have non-identically zero rates or be of the null variety.

The total transition rate will be denoted $r(\eta,x,y)=r_u(\eta,x,y)+r_d(\eta,x,y)$. Note that $r_u$ and $r_d$ are never simultaneously nonzero for any given configuration since we are restricting ourselves to processes which are described by a single interaction map. This is not true in general for multitype contact processes. %See \cite{CoexComp_DurNeu} for an example of a multitype contact process where particle type 3 influences particle type 0 to become type 1 and also to become type 2, thus type 3's influence on type 0 is not unique. This amounts to the need for more than one interaction map to formulate a model with non-unique influences. A sufficiency condition for these types of models to be monotone is given at the end of this paper, but results here do not directly apply.

%{\em Example:}
%For the nearest neighbor contact process, the transition rates are
%\begin{subequations}
%\begin{eqnarray*}
%r_u(\eta,x,y)&=&\lambda_{01}\phi(x,y)I_{\{(0,1)\}}(\eta(x),\eta(y))\\
%r_d(\eta,x,y)&=&\phi(x,y) \big(I_{\{(1,0)\}}(\eta(x),\eta(y))+I_{\{(1,1)\}}(\eta(x),\eta(y))\big)\\
%&=&\phi(x,y)I_{\{1\}}(\eta(x))
%\end{eqnarray*}
%\end{subequations}
%where
%\begin{equation*}
%\phi(x,y) =\left\{
%\begin{aligned}
%&1 &\text{ if } \hspace{2mm} |x-y|=1 \\
%&0 &\text{ otherwise. }
%\end{aligned}\right.
%\end{equation*}
%and $\lambda_{01}=\lambda$. Thus this formulation is equivalent to the standard flip rate formulation of the contact process.

The generator for an IMS is defined on $f\in C(\Xi)$ by
\begin{equation}
G f(\eta)= \displaystyle\sum_{x,y} r(\eta,x,y)(f(\eta^{xy})-f(\eta)).
\label{eq:IMPSgen}
\end{equation}
The state $\eta^{xy}$ represents the state $\eta$ with the particle at site $x$ changed according to the influence of the particle at site $y$:
\begin{equation*}
\eta^{xy}(z) =\left\{
\begin{aligned}
&\eta(z)   &\text{ if } \hspace{2mm} z \neq x \\
&\mathscr{J}(\eta(x),\eta(y)) &\text{ if } \hspace{2mm} z = x
\end{aligned}\right..
\end{equation*}
%These transition rates are assumed to have finite range, $\exists \rho < \infty$ such that $\phi(x,y)=0$ for $|x-y|>\rho$, and be bounded, $\sup_x\sum_{y} r(\eta,x,y)<\infty$. %This holds true since
%$$\sup_x\sum_{y} r(\eta,x,y)=\sup_x\sum_y\lambda_{\eta(x)\eta(y)}\phi(x,y)\le \max_{x,y}\phi(x,y)\cdot \max_{a,b}\lambda_{ab}.$$
%This is obvious now since the $\lambda_{ab}$'s are only numbers and $\phi$ is bounded.

\subsection{The Coupled Rates}
Suppose we have two interaction map systems, $\eta_{\bf t}$ and $\xi_{\bf t}$ on the state space $\Xi$, with transition rates $r_1(\eta,x,y)$ and $r_2(\xi,x,y)$ respectively, with $r_i(\cdot,x,y)=r_{iu}(\cdot,x,y)+r_{id}(\cdot,x,y)$. The coupling described below is essentially the same as the Vasershtein coupling \cite{vaser69} given for spin systems, also known as the basic coupling \cite{LiggettIPS}. The coupled process, $(\eta_{\bf t}, \xi_{\bf t})$, is a Feller process whose state space is $\Xi\times \Xi$ and evolves according to the following rates:
\begin{equation}
(\eta, \xi) \rightarrow
\left\{
\begin{aligned}
& (\eta^{xy}, \xi^{xy}) & \text{ at rate } \hspace{2mm} & \tilde{r}(\eta,\xi,x,y) \\%\label{r3_tilde}\\
& (\eta, \xi^{xy}) & \text{ at rate } \hspace{2mm} & r_2(\xi,x,y)-\tilde{r}(\eta,\xi,x,y) \\%\label{r2_tilde}\\
& (\eta^{xy}, \xi) & \text{ at rate } \hspace{2mm} & r_1(\eta,x,y)-\tilde{r}(\eta,\xi,x,y) %\label{r1_tilde}
\end{aligned}
\right.
\label{eq:coupRates}
\end{equation}
where $\tilde{r}(\eta,\xi,x,y)=\min(r_{1u}(\eta,x,y),r_{2u}(\xi,x,y))+\min(r_{1d}(\eta,x,y),r_{2d}(\xi,x,y))$.
The generator for the coupled process is defined for $f\in C(\Xi\times \Xi)$ by
\begin{subequations}
\begin{eqnarray}
\tilde{G}f(\eta,\xi)= &\displaystyle\sum_{x,y} &(r_1(\eta,x,y)-\tilde{r}(\eta,\xi,x,y))(f(\eta^{xy},\xi)-f(\eta,\xi)) \label{eq:coupGen1}\\
&+\displaystyle\sum_{x,y} &(r_2(\xi,x,y)-\tilde{r}(\eta,\xi,x,y))(f(\eta,\xi^{xy})-f(\eta,\xi)) \label{eq:coupGen2}\\
&+\displaystyle\sum_{x,y} &\tilde{r}(\eta,\xi,x,y)(f(\eta^{xy},\xi^{xy})-f(\eta,\xi)) \label{eq:coupGen3}
\end{eqnarray}
\end{subequations}
This coupling is shown to preserve the partial ordering under certain conditions.

\section{Monotonicity and Attractiveness}
Now we are ready to prove that the coupling preserves the partial ordering on the state space almost surely when the interaction map is attractive and transition rates satisfy certain conditions. This next proof is almost identical to that of Theorem 1.5 in Chapter 3 of \cite{LiggettIPS}, but since the processes here are slightly more complicated a bit more work is necessary.

\begin{theorem} (Extension of Theorem 1.5 in Chapter 3 of \cite{LiggettIPS})
%Define the closed set $K=\{(\eta,\xi)\in \Xi\times \Xi : \eta\le\xi\}$. 
Suppose that the processes $\eta_{\bf t}$ and $\xi_{\bf t}$ on the state space $\Xi$ share the same attractive interaction map. Furthermore suppose that whenever $\eta \le \xi$, the transition rates satisfy
\begin{subequations}
\begin{eqnarray}
& r_{1u}(\eta,x,y) \le r_{2u}(\xi,x,y) &\text{when } \ \eta^{xy}(x)>\xi(x),\label{eq:uprestrpk1} \\
& r_{1d}(\eta,x,y) \ge r_{2d}(\xi,x,y) &\text{when } \ \xi^{xy}(x)<\eta(x)\label{eq:downrestrpk1}.
\end{eqnarray}
\end{subequations}
%Then for all $(\eta,\xi) \in K$ and $t\ge 0$,
Then for all $\eta\le\xi$ and $t\ge 0$,
\begin{equation}
%P^{(\eta,\xi)}[(\eta_t,\xi_t)\in K]=1.
P^{(\eta,\xi)}[\eta_t\le\xi_t]=1.
\end{equation}
\begin{proof} %Define $\mathscr{A}\subset C(\Xi\times \Xi)$ to be the subset of functions which are non-negative everywhere and $0$ on $K$. Since $\mathscr{R}(I-\lambda \tilde{G})=C(\Xi\times \Xi)$ (Proposition 2.8 in Chapter 1 of \cite{LiggettIPS}), there exists $h\in\mathscr{D}(\tilde{G})$ such that $(I-\lambda \tilde{G})h=f$ for each $f\in\mathscr{A}$ and $\lambda\ge0$. Fix $f\in\mathscr{A}$, and define $h$ by $h-\lambda \tilde{G}h=f$. Because $K$ is compact, there exists an $(\eta,\xi)\in K$ where $h$ achieves its maximum. The following cases will show that the generator applied to $h(\eta,\xi)$ is strictly non-positive.
Suppose we have states $\eta\le\xi$, then all that is necessary is to show that the coupled process preserves the partial ordering almost surely. Note that the coupling only allows simultaneous transitions if they are both up or both down.
\paragraph{Case 1} If the lower configuration may jump above the upper, $\xi(x)<\eta^{xy}(x)$, then $r_{1u}(\eta,x,y) \le r_{2u}(\xi,x,y)$ by (\ref{eq:uprestrpk1}). Since $$\tilde{r}(\eta,\xi,x,y)=\min(r_{1u}(\eta,x,y),r_{2u}(\xi,x,y)),$$ we get $r_1(\eta,x,y)-\tilde{r}(\eta,\xi,x,y)=0$. So the problem transition, $(\eta,\xi)\rightarrow(\eta^{xy},\xi)$, occurs at rate zero. %The term (\ref{eq:coupGen1}) disappears from the generator and all that is left are transitions which remain in the set $K$. Because $(\eta,\xi)$ is where $h$ achieves its maximum, the remaining terms of the generator, (\ref{eq:coupGen2}) and (\ref{eq:coupGen3}) are non-positive.
\paragraph{Case 2} If the upper configuration may jump below the lower, $\xi^{xy}(x)<\eta(x)$, then $r_{1d}(\eta,x,y) \ge r_{2d}(\xi,x,y)$ by (\ref{eq:downrestrpk1}). Then $r_2(\eta,x,y)-\tilde{r}(\eta,\xi,x,y)=0$ since $r_2$ is the minimum here leading to $(\eta,\xi)\rightarrow(\eta,\xi^{xy})$ at rate zero. %Once again, similar to the above calculation, the term (\ref{eq:coupGen2}) is now zero. The remaining terms are non-positive since $h$ achieves its maximum at $(\eta,\xi)$.
\paragraph{Case 3} If $\xi(x)<\xi^{xy}(x)$ and $\eta(x)<\eta^{xy}(x)$, then $\eta^{xy}(x)\le\xi^{xy}(x)$ by the first attractive interaction map property. Assuming that $\xi(x)<\eta^{xy}(x)$ would put us back in case 1 and we are done. If $\eta^{xy}(x)\le\xi(x)$, then there are no problem transitions that break the partial ordering. %Once again since $(\eta,\xi)$ is where $h$ attains its maximum, the generator remains non-positive in this case.
\paragraph{Case 4} If $\xi^{xy}(x)<\xi(x)$ and $\eta^{xy}(x)<\eta(x)$, then $\eta^{xy}(x)\le\xi^{xy}(x)$ by the second attractive interaction map property. If $\xi^{xy}(x)<\eta(x)$, then we are back in case 2, otherwise all transitions preserve the partial order.% and the generator is non-positive.
\paragraph{Case 5} If $\eta(x)<\eta^{xy}(x)$ and $\xi^{xy}(x)\le\xi(x)$, then $\eta^{xy}(x)\le\xi(x)$ by the third attractive interaction map property. So $\tilde{r}(\eta,\xi,x,y)=0$ since the minimum up transition rate is $0$ as is the minimum down transition rate. This gives $(\eta,\xi)\rightarrow(\eta^{xy},\xi^{xy})$ at rate zero, and the order is preserved.
\paragraph{Case 6} If $\eta(x)\le\eta^{xy}(x)$ and $\xi^{xy}(x)<\xi(x)$, then $\eta(x)\le\xi^{xy}(x)$ by the fourth attractive interaction map property. So $\tilde{r}(\eta,\xi,x,y)=0$ since the minimum up transition rate is $0$ as is the minimum down transition rate. This gives and $(\eta,\xi)\rightarrow(\eta^{xy},\xi^{xy})$ at rate zero once again preserving the partial ordering.

%This shows that for all possible cases, $\lambda\tilde{G}h(\eta,\xi)\le0$ so that $h(\eta,\xi)\le h(\eta,\xi)-\lambda\tilde{G}h(\eta,\xi)=f(\eta,\xi)=0$. Since $\min_{(\zeta_1,\zeta_2)\in \Xi\times \Xi} h(\zeta_1,\zeta_2) \ge \min_{(\zeta_1,\zeta_2)\in \Xi\times \Xi} f(\zeta_1,\zeta_2)$ (Proposition 2.8 with Definition 2.1 in Chapter 1 of \cite{LiggettIPS}), we see that $h=0$ on $K$, concluding that $h\in\mathscr{A}$. Since $(I-\lambda \tilde{G})^{-1}$ maps $\mathscr{A}$ to itself, so does $\tilde{S}(t)$ (Hille-Yoshida Theorem). Since this is true for any $f\in\mathscr{A}$, when $(\eta,\xi)\in K$, $P^{(\eta,\xi)}[(\eta_t,\xi_t)\in K]=1$ for any $t$.
\end{proof}
\label{thm:PK1}
\end{theorem}

%With the assumptions of Theorem \ref{thm:PK1}, if $\mu_1$ and $\mu_2$ are probability measures on $\Xi$ with $\mu_1\le\mu_2$, then $S_s(t)\mu_1\le S_2(t)\mu_2$ by Corollary 1.7 in Chapter 3 of \cite{LiggettIPS}. This is used to prove when an IMS is monotone.

%Monotonicity for an interaction map system requires restrictions on the interaction map in addition to restrictions on the transitions rates. First we define when an IMS is called attractive.
\begin{definition}
An IMS will be called \textit{attractive} if it has an attractive particle interaction map and given $\eta \le \xi$, the transition rates satisfy:
\begin{subequations}
\begin{eqnarray*}
& r_u(\eta,x,y) \le r_u(\xi,x,y) &\text{when } \ \mathscr{J}(\eta(y),\eta(x))>\xi(x), \label{eq:upRateRest}\\
& r_d(\eta,x,y) \ge r_d(\xi,x,y) &\text{when } \ \mathscr{J}(\xi(y),\xi(x))<\eta(x). \label{eq:downRateRest}
\end{eqnarray*}
\end{subequations}
\label{def:attrRates}
\end{definition}
%Definition \ref{def:attrRates} is similar to that for an attractive spin system given by Liggett in \cite{LiggettIPS}.
Consider two ordered configurations, $\eta\le\xi$. In the event that an interaction may cause $\eta$ to jump above $\xi$ at the site $x$, the latter configuration must have an equal or larger up transition rate and a jump transition that goes at least as far: $\xi(x)<\eta^{xy}(x)\le\xi^{xy}(x)$. When the upper configuration, $\xi$, could possibly jump below the lower configuration, $\eta$, a similar statement applies. A coupling which preserves the partial ordering of the underlying processes is called a monotone coupling. %Now we prove that this definition of attractive is equivalent to being monotone.

\begin{theorem} (Extension of Theorem 2.2 in Chapter 3 of \cite{LiggettIPS}):
An IMS is monotone if and only if it is attractive.
\begin{proof}
%Assuming the process is attractive and given two arbitrary and ordered initial distributions on $\Xi$, $\mu_1\le\mu_2$, we will show that $\mu_1S(t)\le\mu_2S(t)$ which is equivalent to monotonicity of the process by (\ref{eq:monSemi}) and (\ref{eq:monCoup}). For the coupled process (\ref{eq:coupRates}), we are assured the existence of the initial distribution, $\nu$ on $\Xi\times \Xi$, according to Theorem 2.4 in Chapter 2 of \cite{LiggettIPS}. Since $\tilde{S}(t)$ is the semigroup of the coupled process, we know $\nu \tilde{S}(t)(\{(\eta,\xi):\eta\le\xi\})=1$ by Theorem \ref{thm:PK1}. The marginal distributions of $\nu \tilde{S}(t)$ are $\mu_1 S(t)$ and $\mu_2 S(t)$ which proves that $\mu_1 S(t)\le\mu_2 S(t)$ by Theorem 2.4 in Chapter 2 of \cite{LiggettIPS}, completing the proof that the attractive process is monotone.

If the process is attractive, then we have a coupling which preserves the partial order that proves that the IMS is monotone. Assuming that the process is monotone we must prove that it is attractive. We fix $x$ and choose $\eta$ and $\xi$ such that $\eta(y)=\xi(y)$ for all $y\neq x$. Choose $b$ such that $\eta(x)\le\xi(x)<b$. This will be used to prove conditions on up transitions. If no such $b$ exists, there is no problem, as $\eta(x)$ cannot jump above $\xi(x)$. Similarly if there is a $b$ such that $b\le\eta(x)\le\xi(x)$, then this is used to show something about down transitions.

We start with the monotone function $f_{bx}(\eta)=I_{[b..n]}(\eta(x))$, the indicator on all particle values bigger than or equal to $b$. %Note that $b$ is chosen so that $f_{bx}(\eta)=f_{bx}(\xi)$. %This leads to the following calculation.
%The generator and semigroup are related by:
%$$
%G f(\eta) = \sum_x \sum_{y} r(\eta,x,y)(f(\eta^{xy})-f(\eta))=\lim_{t\searrow 0} \frac{S(t)f(\eta)-f(\eta)}{t}
%$$
%For the monotone function $f_{b\tilde{x}}(\eta)$, we don't need to sum over $x$ since it is fixed at $\tilde{x}$ giving:
%$$
%G f_{bx}(\eta) = \sum_{y} r(\eta,x,y)(f_{bx}(\eta^{xy})-f_{bx}(\eta))=\lim_{t\searrow 0} \frac{S(t)f_{bx}(\eta)-f_{bx}(\eta)}{t}
%$$
%Now that $\tilde{x}$ is fixed, we will drop the tildes.
Since $f_{bx}$ is monotone, $S(t)f_{bx}(\eta)\le S(t)f_{bx}(\xi)$. Noting that $b$ was chosen so that $f_{bx}(\eta)=f_{bx}(\xi)$ shows that:
%$$S(t)f_{bx}(\eta)\le S(t)f_{bx}(\xi)$$
%and thus
$$\frac{S(t)f_{bx}(\eta)-f_{bx}(\eta)}{t}\le\frac{S(t)f_{bx}(\xi)-f_{bx}(\xi)}{t}.$$
Taking the limit $t\rightarrow0$ gives:
$$
G f_{bx}(\eta) \le G f_{bx}(\xi).
$$
Plugging in the form of the generator \ref{eq:IMPSgen} results in:
\begin{equation}\begin{aligned}
\sum_{y\in N(x)} r(\eta,x,y)&(f_{bx}(\eta^{xy})-f_{bx}(\eta)) \quad\quad\quad\quad\quad\quad\quad\quad \\
&\le\sum_{y\in N(x)} r(\xi,x,y)(f_{bx}(\xi^{xy})-f_{bx}(\xi))
\end{aligned}\label{eq:rateRel}\end{equation}
This gives the total rate that $\eta(x)$ transitions up into the set $[b .. n]=\{a\in S | b\le a\le n\}$ is less than or equal to the total rate that $\xi(x)$ goes up into the set $[b..n]$, and the total rate that $\eta(x)$ leaves the set $[b..n]$ is greater than or equal to the total rate that $\xi(x)$ does the same.

Then, due to the choice of transition rates, $r(\eta,x,y)=\lambda_{\eta(x),\eta(y)}\phi(x,y)$, (\ref{eq:rateRel}) becomes
\begin{equation}\begin{aligned}
\sum_{y}\lambda_{\eta(x),\eta(y)}\phi(x,y)&(f_{bx}(\eta^{xy})-f_{bx}(\eta)) \quad\quad\quad\quad\quad\quad\quad\quad \\
&\le\sum_{y}\lambda_{\xi(x),\xi(y)}\phi(x,y)(f_{bx}(\xi^{xy})-f_{bx}(\xi))
\end{aligned}\label{eq:rateRel2}\end{equation}
%If we choose a particular neighbor of interest, say $\tilde{y}$ and make $\eta(y)=\eta(\tilde{y})$ and $\xi(y)=\xi(\tilde{y})$ for all $y\neq x$, (\ref{eq:rateRel2}) becomes
%\begin{equation*}
%\lambda_{\eta(x),\eta(\tilde{y})}(x)\sum_{y}\phi(x,y)(f_{bx}(\eta^{x\tilde{y}})-f_{bx}(\eta))\le
%\lambda_{\xi(x),\xi(\tilde{y})}(x)\sum_{y}\phi(x,y)(f_{bx}(\xi^{x\tilde{y}})-f_{bx}(\xi)) \label{eq:rateRel3}
%\end{equation*}
Because the same particle sits at all $y$, $\sum_{y}\phi(x,y)$ cancels out from both sides.
%\begin{equation*}
%\lambda_{\eta(x),\eta(y)}(x)(f_{bx}(\eta^{xy})-f_{bx}(\eta))\le
%\lambda_{\xi(x),\xi(y)}(x)(f_{bx}(\xi^{xy})-f_{bx}(\xi))
%\end{equation*}
%then dropping the tilde on $\tilde{y}$ since this did not depend on exactly which neighbor was singled out:
\begin{equation}\begin{aligned}
\lambda_{\eta(x),\eta(y)}(f_{bx}(\eta^{xy})&-f_{bx}(\eta)) \quad\quad\quad\quad\quad\quad\quad\quad \\
&\le\lambda_{\xi(x),\xi(y)}(f_{bx}(\xi^{xy})-f_{bx}(\xi))
\end{aligned}\label{eq:rateRelnosum}\end{equation}

%\noindent\textbf{Inequalities:}\\
The following inferences come from looking at the possibilities for (\ref{eq:rateRelnosum}).
\begin{packed_enum}
\renewcommand{\labelenumi}{(\roman{enumi})}
\item If $\eta(x)\le\xi(x)<b$, then
\begin{packed_enum}
    \item $\eta^{xy}(x)\ge b \Rightarrow \xi^{xy}(x)\ge b$
    \item $\xi^{xy}(x)< b \Rightarrow \eta^{xy}(x)< b$
\end{packed_enum}
\item If $b\le\eta(x)\le\xi(x)$, then
\begin{packed_enum}
    \item $\xi^{xy}(x)< b \Rightarrow \eta^{xy}(x)< b$
    \item $\eta^{xy}(x)\ge b \Rightarrow \xi^{xy}(x)\ge b$
\end{packed_enum}
\end{packed_enum}
Letting $b=\xi(x)+1$ in inference (i,a), we see that when $\eta^{xy}(x)>\xi(x)$, (\ref{eq:rateRelnosum}) becomes $\lambda_{\eta(x),\eta(y)}\le\lambda_{\xi(x),\xi(y)}$ which is our first rate restriction. For the second rate restriction, let $b=\eta(x)$ in inference (ii,a), we see that when $\eta(x)>\xi^{xy}(x)$, (\ref{eq:rateRelnosum}) becomes $\lambda_{\eta(x),\eta(y)}\ge\lambda_{\xi(x),\xi(y)}$ finishing the rate restrictions for being an attractive process.

If both states have possible up transitions, $\eta(x)<\eta^{xy}(x)$ and $\xi(x)<\xi^{xy}(x)$, then letting $b=\eta^{xy}(x)$ shows that under (i,a), $\eta^{xy}(x)\le\xi^{xy}(x)$. Similarly if both states have possible down transitions, $\eta^{xy}(x)<\eta(x)$ and $\xi^{xy}(x)<\xi(x)$, then $b=\xi^{xy}(x)+1$ under (ii,a) shows that $\eta^{xy}(x)\le\xi^{xy}(x)$. These give the first two requirements for an attractive interaction map.

Assuming $\eta(x)<\eta^{xy}(x)$, $\xi^{xy}(x)\le\xi(x)$, and $b=\xi(x)+1$ along with (i,b) shows that $\eta^{xy}(x)<b$.
If $\eta(x)\le\eta^{xy}(x)$ and $\xi^{xy}(x)<\xi(x)$, then letting $b=\eta(x)$ shows that $b\le\xi^{xy}(x)$ by (ii,b). Now the last two attractive interaction map requirements are met. Since this did not depend on the particular choices of particle types that sat at $x$ and $y$, we are done.
\end{proof}
\label{thm:attrMon}
\end{theorem}

The proof of Theorem \ref{thm:attrMon} is somewhat more involved than the case for spin systems since particle values are allowed to increase by more than one in an interaction giving rise to the possibility of `jump overs'. This result applies to any interaction map system. The inclusion of spatial and temporal variations in $\phi(x,y)$ or in the rate parameters does not affect this result so long as the conditions for being attractive hold over the entire lattice at all times.

\subsection{Reordering the Particles}
If one develops a model which either does not have an attractive interaction map, or the rate restrictions which allow this model to be attractive are not desirable, a re-ordering of the particle values may give an attractive model or more desirable rate restrictions. %In reference to what the model is being used to study, the ordering of the particles may not necessarily have a meaning, it can purely be a mathematical construction for the purpose of creating a monotone model for the benefit of using CFTP.

The basic two-type contact process can be made attractive with a particle re-ordering. The issue that the interaction map is not attractive is resolved by making the permutation $\{0,1,2\}\rightarrow\{1,0,2\}$, in the sense that now $1<0<2$, the interaction map is then attractive. To avoid confusion on such an awkward ordering of integers, we label particle 0 as one species (species 0), particle 1 as empty, and particle 2 remains labeled species 2. The transition rates are still as before: empty sites become species $i$ at rate $\beta_i$ times the number of species $i$ nearby, and species $i$ dies at constant rate $\delta_i$. This information is summarized in Figure \ref{fig:MCPreorder}.

%%%%%%%%%%%%%%%%%%%%%%%%%%%%%
% MCPr Model Transition Rate picture:
\begin{figure}[h]
\begin{center}
\begin{picture}(100,50)
%np3 intmap:
\put(80,8){$\begin{array}{| c || c | c | c | c |}
\hline
_b\hspace{-4pt}\rotatebox{135}{\rule{10pt}{0.5pt}}\hspace{-4pt} ^a  & 0 & 1 & 2\\ \hline\hline
0 & 1 & 0 & 1\\ \hline
1 & 1 & 1 & 1\\ \hline
2 & 1 & 2 & 1\\ \hline
\end{array}$}
%np3 rates:
\put(-50,8){$\begin{array}{| c || c | c | c | c |}
\hline
_b\hspace{-4pt}\rotatebox{135}{\rule{10pt}{0.5pt}}\hspace{-4pt} ^a  & 0 & 1 & 2\\ \hline\hline
0 & \delta_0 & \beta_0 & \delta_2\\ \hline
1 & \delta_0 & \emptyset & \delta_2\\ \hline
2 & \delta_0 & \beta_2 & \delta_2\\ \hline
\end{array}$}
\end{picture}
\end{center}
\vspace{10pt}
\caption{Reordering of the multitype contact process for attractiveness. }
\label{fig:MCPreorder}
\end{figure}
% MCPr Model Int Map picture
%%%%%%%%%%%%%%%%%%%%%%%%%%%%%%%
This model now has an attractive interaction map and no rate restrictions. Normally 0 represents an empty site, but this shows that thinking more carefully about particle labels is beneficial. This gives us ordered extremal stationary distributions: $\nu_0\le\delta^{\bar{1}}\le\nu_2$ where $\delta^{\bar{i}}$ is point mass on the state $\eta(x)=i$ for all $x$ and $\nu_i=\displaystyle\lim_{t\rightarrow\infty}\delta^{\bar{i}}S(t)$ is the invariant measure for species $i$ in isolation.

%This is one problem with using integers to represent the particle types; it gives us a preconceived notion on how they should be ordered. This method was discovered by studying a model where the particles were not labeled with integers, thus no natural ordering was assumed.

\section{The Graphical Representation}
Now we define a graphical coupling which couples all states simultaneously. This graphical representation was first introduced by Harris \cite{Harris72,Harris78} and subsequently modified by others \cite{FristedtGray,LiggettIPS,Liggett99}. This method is useful for exact sampling algorithms such as CFTP. THe graphical coupling here is built from that in \cite{FristedtGray}. %useful since it allows multiple copies of the process to be constructed simultaneously. %The construction here is derived mostly from that given in chapter 32 of Fristedt and Gray\cite{FristedtGray}.

For each site, $x$, and every $y$ such that $\phi(x,y)>0$, let $U_{xy}$ and $D_{xy}$ be two independent, identically distributed Poisson point processes on $(0,\infty)\times(0,\infty)$ with intensity equal to two-dimensional Lebesgue measure. Assume that our rates have an upper bound, $c$. For each $x$, define $\mathfrak{T}_x=\{T_{x,1}< T_{x,2}< ...\}$ by $T_{x,0}=0$ and
$$T_{x,n}=\inf\{t>T_{x,n-1} : (v,t)\in \bigcup_y U_{xy}\cup D_{xy} \text{ for some } v \le c\}$$
which will represent the times at which transitions could possibly occur at the site $x$. This is a projection of a union of independent Poisson point processes, and the intensity measure of points in $\mathfrak{T}_x$ is $2c\#N(x,\rho)$ where $N(x,\rho)=\{y:|x-y|\le\rho \text{ and } \phi(x,y)>0\}$, so $\mathfrak{T}_x$ is also a Poisson point process.

\begin{figure}[h]
\begin{center}
\begin{picture}(200,220)
\put(-30,100){\vector(0,1){50}}
\put(-40,90){time}
\put(100,-24){$x$}
\put(-8,-14){$-2$}
\put(33,-14){$-1$}
\put(78,-14){$0$}
\put(118,-14){$1$}
\put(158,-14){$2$}
\put(198,-14){$3$}
\put(0,0){\rule{200pt}{2pt}}
\put(0,0){\rule{.5pt}{3in}}
\put(40,0){\rule{.5pt}{3in}}
\put(80,0){\rule{.5pt}{3in}}
\put(120,0){\rule{.5pt}{3in}}
\put(160,0){\rule{.5pt}{3in}}
\put(200,0){\rule{.5pt}{3in}}
\put(0,60){\vector(1,0){36}}\put(40,60){\circle*{8}}\put(46,60){$1.2$}
\put(40,130){\vector(-1,0){36}}\put(0,130){\circle{8}}\put(-20,130){$0.4$}
\put(40,20){\vector(1,0){36}}\put(80,20){\circle*{8}}\put(86,20){$0.7$}
%\put(40,170){\vector(1,0){36}}\put(80,170){\circle{8}}
%
\put(120,100){\vector(-1,0){36}}\put(80,100){\circle{8}}\put(60,100){$1.8$}
\put(80,190){\vector(1,0){36}}\put(120,190){\circle*{8}}\put(126,190){$2.3$}
\put(160,70){\vector(-1,0){36}}\put(120,70){\circle{8}}\put(100,70){$0.1$}
\put(120,110){\vector(1,0){36}}\put(160,110){\circle*{8}}\put(166,110){$1.9$}
\put(160,150){\vector(-1,0){36}}\put(120,150){\circle{8}}\put(100,150){$2.0$}
\end{picture}
\end{center}
\vspace{10pt}
\caption{A realization of the graphical representation of the point process coupling. All points with $u\geq 3$ have been left out. Closed and open circles represent possible down and up transition points respectively.}
\label{fig:graphrep}
\end{figure}
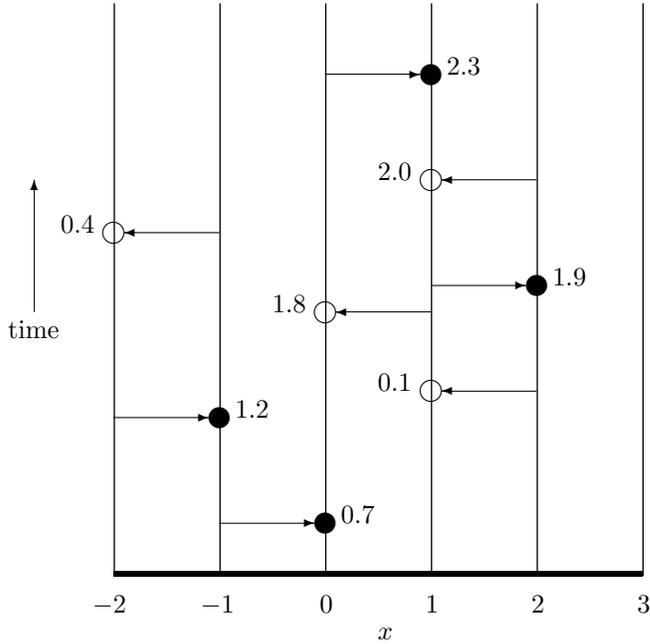

The graphical representation is created on a grid $\Lambda \times [0,\infty)$. At each point, $t$, in $\mathfrak{T}_x$, if $\exists u, y$ such that $(u,t) \in U_{xy}$ draw an arrow from $y$ pointing to $x$ with an open circle at $x$. If $\exists u, y$ such that $(u,t) \in D_{xy}$ draw an arrow from $y$ pointing to $x$ with a closed circle at $x$. Write the $u$ values next to the tip of each arrow. Figure \ref{fig:graphrep} give a realization of the graphical representation for a one dimensional index set and rates bound above by $c=3$.
%%%%%%%%%%%%%%%%%%%%%%%%%%%%%
% graphical coupling picture:
% end graphical coupling picture
%%%%%%%%%%%%%%%%%%%%%%%%%%%%%%%

These point processes are used to evolve the interaction map system. Given an initial state, $\eta$, the following theorem describes how the graphical coupling evolves the process over time.
\begin{theorem}The path $\eta_\bft$ is constructed according to the following rules is the interaction map system with interaction map $\mapJ$, and generator given by (\ref{eq:IMPSgen}).
\begin{packed_enum}
\item Up transition rule: The particle at site, $x$, is replaced by the particle given by
$\mathscr{J}(\eta_{t-}(x),\eta_{t-}(y))$
at time $t\in\mathfrak{T}_x$ if there exists a $u$ such that $(u,t)\in U_{xy}$, with
$\mathscr{J}(\eta_{t-}(x),\eta_{t-}(y))>\eta_{t-}(x),$
and $u\le r_u(\eta_{t-},x,y)$.
\item Down transition rule: The particle at site, $x$, is replaced by the particle given by
$\mathscr{J}(\eta_{t-}(x),\eta_{t-}(y))$
at time $t\in\mathfrak{T}_x$, if there exists a $u$ such that $(u,t)\in D_{xy}$, with
$\mathscr{J}(\eta_{t-}(x),\eta_{t-}(y))<\eta_{t-}(x),$
and $u\le r_d(\eta_{t-},x,y)$.
\end{packed_enum}
\label{thm:UC}
\end{theorem}
The construction of these Poisson point processes and Theorem \ref{thm:UC} is based upon the construction in Chapter 32 of \cite{FristedtGray}. This is the basis for an accurate method of simulating these types of processes.

\begin{theorem}
The graphical construction in Theorem \ref{thm:UC} maintains the partial order of an attractive IMS.
\label{thm:graphCoupMon}
\end{theorem}
Theorem \ref{thm:graphCoupMon} is proven by the fact that all transitions preserve the order of configurations. %The only possibility for the partial ordering to be possibly broken is if up and down transitions are allowed simultaneously, this is not the case here due to the segregation of the up and down point processes. For any attractive interaction map system, it is easy to construct a monotone coupling using this method. There are other graphical couplings as well. Transitions can be grouped into point processes according to any criteria, they just should not be grouped together in such a way that allows an up transition for one configuration and a down transition for another if these two transitions break the ordering.

\begin{theorem}
Consider two sets of parameters satisfying $\lambda^{(1)}_{ab}\le\lambda^{(2)}_{ab}$ for all up transitions and $\lambda^{(1)}_{ab}\ge\lambda^{(2)}_{ab}$ for all down transitions. If two states satisfying $\xi^{(1)}\le\xi^{(2)}$ are the initial states for the processes with the corresponding parameter sets above, then $\xi^{(1)}_t\le\xi^{(2)}_t$ for all $t\ge0$ for the above graphical construction.
\label{thm:paramMon}
\end{theorem}
Theorem \ref{thm:paramMon} is proven by applying Theorem \ref{thm:PK1} to the graphical coupling. This shows that an attractive interaction map process is monotone in each of its rate parameters. This is similar to the monotonicity theorems in \cite{CoexComp_DurNeu} and \cite{mScp_Krone}.

\section{Multiple Interaction Maps}
One may desire to include more complicated interactions such as non-unique interaction, i.e. when a particular particle type can have multiple influences over another, or if a constant transition rate is included along with several other particle interactions. Take for example the grass--bushes--trees model studied in \cite{DurSwin91}.

In order to make this model attractive, we reorder the particles as for the basic two-type contact process and introduce a separate interaction map for the births of trees. The basic interaction map, $\mapJ_1$, and transition rates are given by Figure \ref{fig:GBT01}. This interaction map is the basic two-type contact process map with births for species two removed and is attractive with no rate restrictions. Now we just need to account for this extra birth event. This amounts to including the extra interaction map, $\mapJ_2$, given by Figure \ref{fig:GBT02}, and this map is attractive as well.
%%%%%%%%%%%%%%%%%%%%%%%%%%%%%
% MCPr Model Transition Rate picture:
\begin{figure}[h]
\begin{center}
\begin{picture}(100,50)
%np3 intmap:
\put(80,8){$\begin{array}{| c || c | c | c | c |}
\hline
_b\hspace{-4pt}\rotatebox{135}{\rule{10pt}{0.5pt}}\hspace{-4pt} ^a  & 0 & 1 & 2\\ \hline\hline
0 & 1 & 0 & 1\\ \hline
1 & 1 & 1 & 1\\ \hline
2 & 1 & 1 & 1\\ \hline
\end{array}$}
%np3 rates:
\put(-50,8){$\begin{array}{| c || c | c | c | c |}
\hline
_b\hspace{-4pt}\rotatebox{135}{\rule{10pt}{0.5pt}}\hspace{-4pt} ^a  & 0 & 1 & 2\\ \hline\hline
0 & \delta_1 & \beta_1 & \delta_2\\ \hline
1 & \delta_1 & \emptyset & \delta_2\\ \hline
2 & \delta_1 & \emptyset & \delta_2\\ \hline
\end{array}$}
\end{picture}
\end{center}
\vspace{10pt}
\caption{Grass--bushes--trees initial rate parameters and interaction map, $\mapJ_1$.}
\label{fig:GBT01}
\end{figure}
% MCPr Model Int Map picture
%%%%%%%%%%%%%%%%%%%%%%%%%%%%%%%
%%%%%%%%%%%%%%%%%%%%%%%%%%%%%
% GBT Model extra Int/rate Map picture:
\begin{figure}[h]
\begin{center}
\begin{picture}(100,50)
%np3 intmap:
\put(80,8){$\begin{array}{| c || c | c | c | c |}
\hline
_b\hspace{-4pt}\rotatebox{135}{\rule{10pt}{0.5pt}}\hspace{-4pt} ^a  & 0 & 1 & 2\\ \hline\hline
0 & 0 & 1 & 2\\ \hline
1 & 0 & 1 & 2\\ \hline
2 & 2 & 2 & 2\\ \hline
\end{array}$}
%np3 rates:
\put(-50,8){$\begin{array}{| c || c | c | c | c |}
\hline
_b\hspace{-4pt}\rotatebox{135}{\rule{10pt}{0.5pt}}\hspace{-4pt} ^a  & 0 & 1 & 2\\ \hline\hline
0 & \emptyset & \emptyset & \emptyset\\ \hline
1 & \emptyset & \emptyset & \emptyset\\ \hline
2 & \beta_2 & \beta_2 & \emptyset\\ \hline
\end{array}$}
\end{picture}
\end{center}
\vspace{10pt}
\caption{Extra rate parameters and interaction map, $\mapJ_2$, for the birth of trees.}% onto sites occupied by bushes.}
\label{fig:GBT02}
\end{figure}
% GBT Model extra Int Map picture
%%%%%%%%%%%%%%%%%%%%%%%%%%%%%%%
%There are only two actual transitions occurring in this map, the births of species 2.
This is not a unique formulation, deaths could be distributed among  both interaction maps and still maintain attractiveness. When simulating this process we need three Poisson Point Processes: $U$ = up transitions for $\mapJ_1$ which includes deaths for species zero, $D$ = down transitions for $\mapJ_1$ which includes deaths of species two and births of species zero, and $B$ = birth events for species two (up transitions for $\mapJ_2$). This shows that the grass--bushes--trees model is attractive with no rate restrictions.

The idea is that any number of interaction maps can be used. If each map is attractive, then a collection of rate restrictions allows the model to be monotone. Each map is assigned distinct up and down point processes associated with it in the graphical coupling.
%
%Formulating a process with multiple interaction maps my actually relax some parameter restrictions for monotonicity.
The equivalency of attractiveness and monotonicity for a process not described by a single interaction map is not discussed here in detail, but a sufficiency condition is given.

\begin{theorem}{\bf Monotonicity Sufficiency Condition}:
If an interacting particle system is formulated with multiple attractive interaction maps, $\mathscr{J}_i$ $i=1, ...,n$, and the corresponding attractive transition rates, then $\{U_i, D_i\}_{i=1,...,n}$ is a monotone coupling for the process where $U_i$ and $D_i$ are the point processes for the up and down transitions for interaction map $\mapJ_i$ respectively.
\label{thm:multiMap01}
\end{theorem}
Each of the point processes $\{U_1, D_1, U_2, D_2, ..., U_n, D_n\}$ preserves the partial ordering of the state space, thus we have a monotone coupling. This allows non-unique interactions between two particle types. If multiple interaction maps are not used, monotonicity is given by inequalities involving sums of rate parameters rather than a comparison of individual parameters. This may allow the relaxation of the requirement that each individual map and its parameters be attractive, but this is not pursued further here. The grass--bushes--trees model in Figures \ref{fig:GBT01} and \ref{fig:GBT02} is attractive with no rate parameter restrictions according to this theorem. The two-stage contact process can also be seen to be monotone with no rate restrictions.

With the interaction map formulation presented here, monotone properties of multitype contact processes can be assessed quickly. While the main result only applies to processes with a single interaction map, it is still useful for determining whether or not a process with multiple interaction maps is monotone.

\section*{Acknowledgements}
This work is part of my doctoral dissertation at The University of Arizona in The Program in Applied Mathematics under the supervision of Joseph C. Watkins. %I want to acknowledge his encouragement and dedication to my success during this time and for his invaluable suggestions on how best to complete this project.

\end{document}